\documentclass[12pt]{amsart}

\usepackage{amssymb,amsmath,amsthm,longtable,multirow,array}

\usepackage{datetime}
\usepackage{amsfonts}
\usepackage{color}
\usepackage{mathrsfs}
\newtheorem{Theorem}{Theorem}[section]
\newtheorem{Proposition}[Theorem]{Proposition}

\newtheorem{Lemma}[Theorem]{Lemma}
\newtheorem{Corollary}[Theorem]{Corollary}
\theoremstyle{definition}
\newtheorem{Definition}[Theorem]{Definition}
\theoremstyle{remark}
\newtheorem{Example}[Theorem]{Example}

\newtheorem*{remark*}{Remark}

\providecommand\ba[1]{\begin{align*}#1\end{align*}}

\providecommand\CC{{\mathbb C}}
\providecommand\RR{{\mathbb R}}

\providecommand\QQ{{\mathbb Q}}
\providecommand\NN{{\mathbb N}}

\providecommand\brs{\begin{remark*}}
\providecommand\ers{\end{remark*}}

\providecommand\Om{\Omega}

\providecommand\shc{ $\sigma$-holomorphically-convex\ }
\providecommand\shct{ $\sigma$-holomorphically-convex}

\providecommand\dist{\operatorname{dist}}

\providecommand\be{\begin{enumerate}}
\providecommand\ee{\end{enumerate}}
\providecommand\bT{\begin{Theorem}}
\providecommand\eT{\end{Theorem}}
\providecommand\bP{\begin{Proposition}}
\providecommand\eP{\end{Proposition}}
\providecommand\bD{\begin{Definition}}
\providecommand\eD{\end{Definition}}
\providecommand\bE{\begin{Example}}
\providecommand\eE{\end{Example}}
\providecommand\bL{\begin{Lemma}}
\providecommand\eL{\end{Lemma}}
\providecommand\bC{\begin{Corollary}}
\providecommand\eC{\end{Corollary}}

\providecommand\bpp{\begin{proof}}
\providecommand\epp{\end{proof}}
\providecommand\bee{\begin{equation}}
\providecommand\eee{\end{equation}}

\providecommand\Conv{\operatorname{Conv}}

\providecommand\beqq{\begin{eqnarray*}}
\providecommand\eeqq{\end{eqnarray*}}

\begin{document}
\title[Convergence Sets of Power Series]{On Convergence Sets of Power Series with Holomorphic Coefficients}
\author{Basma Al-Shutnawi}
\author{Hua Liu}
\author{Daowei Ma}

\begin{abstract} We consider convergence sets of formal power
series of the form $f(z,t)=\sum_{n=0}^{\infty} f_n(z)t^n$, where
$f_n(z)$ are holomorphic functions on a domain $\Omega$ in $\CC$. A
subset $E$ of  $\Omega$ is said to be a convergence set in $\Om$ if
there is a series $f(z,t)$ such that $E$ is exactly the set of
points $z$ for which $f(z,t)$ converges as a power series in a
single variable $t$ in some neighborhood of the origin. A
$\sigma$-convex set  is defined to be the union of a countable
collection of polynomially convex compact subsets. We prove that a
subset of $\CC$ is a convergence set if and only if it is
$\sigma$-convex.
\end{abstract}

\keywords{formal power series, analytic functions, convergence sets}
\subjclass[2000]{Primary: 32A05, 30C85}
\address{salmashut@yahoo.com,  Department of Mathematics, Tafila technical university,\newline P.O.~Box~179, Tafila 66110, Jordan}
\address{hualiu@tute.edu.cn, Department of Mathematics, Tianjin
University of Technology and Education, Tianjin 300222, China}
\address{dma@math.wichita.edu, Department of Mathematics, Wichita
State University, Wichita, KS 67260-0033, USA}

\maketitle

\bigskip
\section{Introduction}
The purpose of this article is to  describe  the convergence sets of formal power
series with holomorphic coefficients. The study of convergence sets comes from  generalizations of Hartogs Theorem (see \cite{Lelong 1951,AM,S,Levenberg and Molzon 1988, spallek}).
Our approach is motivated by recent work  \cite{B.M,B.D.N,Ma T} concerning formal power series $F(z,t)=\sum_{n=0}^{\infty} P_n(z)t^n$, whose coefficients are polynomials  of one or more complex variables. In some of these studies the authors focus
on  power series of the form
$F(z,t)=f(tz_1,\dots,tz_N)=\sum_{n=0}^{\infty} P_n(z)t^n$, where
 $P_n(z)$ are homogeneous polynomials of degree $n$ for $n\in \NN$. We say that $E\subset\CC^N$ is the convergence set of $F$ if for every $z\in E $ there exists some $r_z>0$ such that $F(z,t)$ converges for $t<r_z$  while for each $z\in  \CC^N\setminus E$ the radius of convergence of $F(z,t)$ equals 0.
Note that it is always assumed that $\deg P_n\leq n$  in these investigations.

Some related problems from operator calculus
and renormalization of quantum field theory ( see, {\it e.g.,} \cite{knill}) are concerned with the formal series of the form
$$ F(T,t)=\sum_{n=0}^\infty K_n(T)t^n,$$
where $K_n(T)$ belong to the $C^\ast$ algebra generated by some
operator $T$. It is necessary to discuss the convergence set of $F(z,t)$ for both spectrum analysis of  $F(T,t)$ and
perturbation theory. Here
$K_n(T)$ are  holomorphic functions on some
neighborhood of the compact set $K$,  the spectra of $T$. It is desirable to find the necessary and sufficient
conditions for   a set $E$ to be the convergence set of some $F(z,t)=\sum_{n=0}^\infty f_n(z)t^n$. This article answers the question completely when
$N=1$.

\section{Convergence Sets}
Let $\Omega$ be an open subset in the complex space $\CC^N$. Denote
by { $\mathscr O(\Om)$} the set of  holomorphic functions on
 $\Omega$.

We consider the power series of the form
\begin{equation}\label{*}
    f(z,t)=\sum^{\infty}_{n=0}f_n(z)t^n,\quad z\in \Omega,
\end{equation}
where $f_n(z)\in\mathscr{O}(\Om)$ and $t$ is a complex variable. We
denote by { $\mathscr{O}(\Om)[[t]]$} the collection of the series of
form (1).

\bD Let $f(z,t)\in\mathscr{O}(\Om)[[t]]$. We define the { {\it
convergence set of $f$ in} $\Om$} by
$$ \Conv_\Om (f)=\{z\in \Omega:\ f(z,t)\  \hbox{converges in some neighborhood
of}\ 0 \}, $$
or equivalently,
$$ \Conv_\Om (f)=\{z\in \Omega:\ |f_n(z)|<r_z^n\  \hbox{for some}\
r_z>0\  \hbox{and every}\ n\in\NN\}. $$ \eD

\bD A subset $E\subset\Omega$ is said to be a {\it convergence set
 in $\Om$} if there exists  an $f\in\mathscr{O}(\Om)[[t]]$ such
that $E=\Conv_\Om(f)$. A convergence set in $\CC^N$ is also simply called  as {\it convergence set}.
\eD

\bP Let $K$ be a polynomially convex compact set in $\CC^N$. Then
$K$ is a convergence set { in $\CC^N$}.\eP

\bpp Let $m$ be any positive integer  and  $y\in \CC^N \setminus K$.
Since $K$ be a polynomially convex,  there exists a polynomial $P_y(z)$ such that $|P_y(y)|>m$ and
$|P_y(z)| \leq 1$ for $z\in K$.

Set $U_{y}=\{x\in\CC^N:|P_{y}(x)|> m\}.$  The open cover $U_{y},
y\in\CC^N\setminus K$, of the set $\CC^N\setminus K$ contains a
countable subcover $ U_{y_{k}},  k=1,2,\ldots$. Now denote by
$P_{mk}(z)=P_{y_{k}}(z)$. For each $m$ we get a sequence
$\{P_{mk}\}_{k=1}^\infty.$
Since the set $\{P_{mk}\}$ is countable we can arrange it as a
sequence $\{h_j(z)\}_{j=1}^ \infty$. Set
$$
f(z,t)=\sum_{j=1}^\infty  h_j^j(z)t^j.
$$

Suppose that $z\in K$. Then for each $j$, $|h_j(z)|\le 1$. Hence
$z\in \Conv_{\CC^N}(f)$. Consequently, $K\subset\Conv_{\CC^N}(f)$.

Now suppose that $z\in \CC^N\setminus K$. Then for each $m\in \NN$
there is a $k\in \NN$ such that $|P_{mk}(z)|\ge m$. It follows that
the sequence $\{|h_\ell(z)|\}$ is unbounded.  So the formal power
series $f(z,t)$ is divergent at  $z$. Consequently,
$\Conv_{\CC^N}(f)\subset K$. That is,
$K$ is a convergence set in $\CC^N$. \epp

The following proposition provides a necessary condition for a set to be a convergence set.

\bP\label{226} Let $E$ be a convergence set in $\Om\subset \CC^N$. Then $E$ is  an $F_\sigma$ set.\eP

 \begin{proof}
Suppose that $E=\Conv_\Om(f)$, where  \ba{
f(z,t)=f_0(z)+f_1(z)t+\cdots+f_n(z)t^n+\cdots\in \mathscr
O(\Om)[[t]].}

For $j\in\NN$ denote by $K_j=\{z\in\Omega:\ \hbox{dist}\{z,\partial \Omega\}\ge \frac 1j\ \hbox{and}\ |z|<j$. Then    $\Om=\cup_{j=1}^\infty K_j$ and each $K_j$ is contained in the interior of $K_{j+1}$. We now prove that
\begin{equation}\label{trial}
E=\bigcup_{j=1}^{\infty}\bigcap_{n=0}^\infty
\{z\in K_j: \;|f_n(z)|\leq j^n\}.
\end{equation}
For $z\in E$, suppose that $z\in K_L$ for some integer $L$. By the definition of
convergence set, there exists a positive integer $J$ such that
$$
|f_n(z)|<J^n,\;n\in\NN.
$$
Let $m=\max\{L,J\}$. Then  $z\in \bigcap_{n=0}^\infty \{z\in K_m: \;|f_n(z)|\leq
m^n\}$. On the other hand, assume that
$z\in
\bigcup_{j=1}^{\infty}\bigcup_{l=1}^{\infty}
\{z\in K_j: \;|f_n(z)|\leq j^n\}$. Then there exist a
positive integer $j$ such that $z\in \bigcap_{n=0}^\infty \{z\in
K_l: \;|f_n(z)|\leq j^n\}$. So $|f_n(z)|\leq j^n$ for all $n$, i.e.,
$z\in E$.

  It is clear that
$\bigcap_{n=0}^\infty \{z\in K_j: \;|f_n(z)|\leq j^n\}$ are closed.
By (\ref{trial}), $E$ is an $F_\sigma$ set.
\end{proof}

The converse of Theorem~\ref{226} is not true, for which we will
give a counterexample in the next section.
\vskip2mm
We now discuss the intersection of several convergence sets.

{ \bP\label{intersection} Let $E_1,\dots, E_k$ be convergence sets in $\Om$. Then  the intersection $E:=\cap_{j=1}^k E_j$ is also a convergence set in $\Om$. \eP}

\begin{proof} It suffices to prove that the intersection of two convergence sets in $\Om$ is a convergence set  in $\Om$.
Suppose that we have two formal power series { in
$\mathscr{O}(\Om)[[t]]$}
 \ba{
f(z,t)=f_0(z)+f_1(z)t+\cdots+f_n(z)t^n+\cdots,} and \ba{
g(z,t)=g_0(z)+g_1(z)t+\cdots+g_n(z)t^n+\cdots} and the corresponding convergence sets
$A=\hbox{Conv}_{\Om}(f)$ and $B=\hbox{Conv}_{\Om}(g)$, respectively.
Define by
\ba{F(z,t)&=F_0(z)+F_1(z)t+\cdots= f(z,t^2)+ tg(z,t^2)\\
            & =f_0(z)+g_0(z)t+f_1(z)t^2+g_1(z)t^3+\cdots.}
Then
\begin{equation}\label{}
  F_n(z)=\left\{\begin{array}{ll}
  f_{\frac n2}(z),&\ n\ \hbox{is even};\\
  g_{\frac{n-1}2}(z),& \ n\ \hbox{is odd}.\end{array}\right.
\end{equation}
For $z\in A\cap B$ suppose that $|f_n(z)|<r_A^n$ and $|g_n(z)|<r_B^n$ for some $r_A,r_B>0$ and every $n\in \NN$. Then $|F_n(z)|<((\max\{r_A,r_B\})^{\frac12})^n$. So $A\cap B\subset \Conv_\Om(F)$.

On the other hand, for $z\in \Conv_\Om(F)$ suppose that $|F_n(z)|<r^n$. Then both $|f_n(z)|$ and $|g_n(z)|$ are less than $(r^2)^n$, i.e., $\Conv_\Om(F)\subset A\cap B.$
\end{proof}

In the rest of the article we only consider the case $N=1$. We
discuss some properties for the convergence sets in the complex
plane.  For any countable set $E$ in $\CC$ we in the following theorem construct a formal
series whose convergence set in $\CC$ is exactly $E$. For the convenience of the proof, denote by  {\it i.e.}, \ba{N_r(S)=\{z\in\CC:
|z-p|<r\;\text{ for some $p\in S$}\}}
for $S\subset \CC$ and $r>0$.


\bT\label{cont conv} Let $S=\{z_1,z_2,\dots\}$ be a countable infinite subset of $\CC$. Define
an $F\in\CC[z][[t]]$ by
$$F(z,t)=\sum_{n=0}^\infty C_{n} \Big[\prod_{j=1}^n (z-z_{j})\Big]t^n, $$
 where $C_{n}={(n/\gamma_{n}})^n$, and
$$\gamma_{n}=\min(\frac{1}{2}\min_{1\leq i<j\leq n+1} |z_{i}-z_{j}|,{1}/{n}).$$
{ Then $\Conv_\CC(F)=S$.}\eT

\begin{proof}
Note that $\gamma_n$ is positive since $z_i$ are pairwise distinct.
Let $L_j=\{z_1,\cdots,z_j\}$ for $j\in\NN$. We
now prove that \bee\label{a1}\bigcap_{j=k}^\infty N_{{\gamma_{j}}}(L_j)=L_k,\ k\in\NN.\eee

 We only
need to prove  $\bigcap_{j=k}^\infty N_{{\gamma_{j}}}(L_j) \subset L_k$ , which would follow from the
following statement:
 \bee\label{a2}\bigcap_{s=k}^\infty N_{{\gamma_{s}}}(L_s)\subset N_{\gamma_j}(L_k),\  \mbox{for} \  j\geq k.\eee
We prove (\ref{a2}) by induction on $j$. It is obvious for $j=k$
since $N_{{\gamma_{k}}}(L_k)=U_k$.
 Suppose the statement is true for $j=N\geq k$. Let
$z\in \bigcap_{s=k}^\infty N_{{\gamma_{s}}}(L_s).$  For $i\not=j$, $1\le i,j\le N+1$,
since $|z_j-z_i|\geq 2\gamma_N$, we get that $N_{\gamma_N}(z_i)\cap
N_{\gamma_N}(z_j)=\emptyset$. It follows that \bee\label{a3}
\Big(\bigcup_{\ell=k+1}^{N+1} N_{\gamma_{N}}(z_{\ell}) \Big) \cap
\Big(N_{\gamma_N}(L_k)\Big)= \emptyset.\eee

By the induction hypothesis $z\in  N_{\gamma_N}(L_k)$ and by
(\ref{a3}), $z\notin \bigcup_{\ell=k+1}^{N+1}
N_{\gamma_{N}}(z_{\ell})$. Since $\gamma_{N+1}\leq \gamma_{N}$ we
see that $\bigcup_{\ell=k+1}^{N+1} N_{\gamma_{N+1}}(z_{\ell})\subset
\bigcup_{\ell=k+1}^{N+1} N_{\gamma_{N}}(z_{\ell})$, hence \bee
\label{a54}z\notin\bigcup_{\ell=k+1}^{N+1}
N_{\gamma_{N+1}}(z_{\ell}).\eee
On the other hand, we know that
\bee\label{a55} z\in N_{{\gamma_{N+1}}}(L_{N+1})=\Big(\bigcup_{\ell=k+1}^{N+1}
N_{\gamma_{N+1}}(z_{\ell}) \Big) \cup
\Big(N_{\gamma_{N+1}}(L_k)\Big).\eee By (\ref{a54}) and (\ref{a55}),
$z\in N_{\gamma_{N+1}}(L_k)$. This completes the induction step, and
therefore the statement is proved.

Now let $P_0(z)=1$ and, for $n\geq1$,
 \begin{equation}
 P_n(z)= C_n \prod_{k=1}^n (z-z_{k}).
 \end{equation}
Then $ F(z,t)=\sum_{n=0}^\infty P_n(z)t^n $.
For $n\geq k$, we have $P_n(z_k)=0$. It follows that
$S\subset\Conv_\CC(F)$.

Now suppose that $z\not\in S$. By (\ref{a1}),
$z\notin\bigcap_{j=k}^\infty U_j$ for each $k$. It follows that
there is a strictly increasing sequence $j_k$ of positive integers
such that
 $z\notin U_{j_k} $ for $ k=1,2,\cdots$. Then
 \begin{equation}
 |P_{j_k}(z)|\geq (\frac{j_k}{\gamma_{j_k}})^{j_k}\prod_{i=1}^{j_k}({\gamma_{i}})\geq({j_{k}})^{j_k},
 \end{equation}
since $\gamma_i\geq \gamma_{j_k}$ for $i \leq j_k.$ This implies
that $z\notin \Conv_\CC(f).$ Therefore $S=\Conv_\CC(f).$
\end{proof}

\bE By Theorem~\ref{cont conv} the set $\mathbb Q$ of rational numbers
is a convergence set in $\CC$. But by Theorem~\ref{226} the set
$\RR\setminus \QQ$ of irrational numbers is not so, since it is not
an $F_{\sigma}$ set. \eE

{ \bT Let $\Omega \subset \CC$, let $ S$ be a countable  dense
subset of $\Omega$, let $\{C_{n}\}$ be a sequence of positive
numbers, and  let { $A=\{a_1,a_2,\cdots,a_k\}$ be any finite subset
of $\Omega$. Then there exists an enumeration $\{z_1, z_2, \dots\}$
of $S$, such that $A\subset \Conv_\Om(F)$, }where $F$ is defined by
$$F(z,t)=\sum_{n=0}^\infty C_{n}\Big[ \prod_{j=1}^n
(z-z_{j})\Big]t^n.$$\eT}

\begin{proof}Suppose that the diameter of $A$ is $d$.  Let $S=\{s_1,s_2,\dots\}$. We choose distinct points $z_1,z_2,\cdots,$  from $S$ such that the following are satisfied:
\begin{enumerate}
  \item $|z_{l(k+1)+i}-a_i|<d,$ and
 $ C_{l(k+1)+i+p}|z_{l(k+1)+i}-a_i|<\frac d{(l+2)!}$, for
 $l=0,1,2,\cdots,\ i=1,2,\cdots,k,\ p=0,1,\cdots,k$;
 \item for $l=1,2,\cdots,\ z_{l(k+1)}=s_{\tau(l)}$, where
 $$\tau(l)=\min\{p\in\NN:|s_p-a_i|<ld,s_p\in S\setminus \{z_1,\cdots,z_{l(k+1)-1}\}\}.$$
\end{enumerate}

We now show that
\begin{equation}
\nonumber |C_n(a_i-z_1)\cdots(a_i-z_n)|<(2d)^n,\ n\geq k+1,\
i=1,\cdots,k.
\end{equation}
\indent Fix $n\geq k+1$ and $1\leq i\leq k$. Choose $l\geq 0$ and $0\leq p\leq k$ such that $n=l(k+1)+i+p$.\\
Since
$$
|a_i-z_{m(k+1)+j}|\leq |a_i-a_j|+|a_j-z_{m(k+1)+j}|<d+d=2d,$$
for $m=1,2,\cdots$ and $1\leq j\leq k$, and since
$$
|a_i-z_{m(k+1)}|<(m+1)d<(m+1)(2d),$$
for $m=1,2,\cdots$, we have




\begin{eqnarray}\label{a}
&\ &|C_n(a_i-z_1)\cdots(a_i-z_n)|\\
 \nonumber&=&|(a_i-z_1)\cdots (a_i-z_{l(k+1)+i-1})|\cdot|C_{l(k+1)+i+p}(a_i-z_{l(k+1)+i})|\cdot  \\
 \nonumber &\ &|a_i- z_{l(k+1)+i+1}|\cdots|a_i-z_{l(k+1)+i+p}| \\
  \nonumber &\leq& (l+2)!(2d)^{l(k+1)+i-1}[\frac d {(l+2)!}](2d)^{p}\\
 \nonumber&=&(2d)^{n-1}(l+2)!\frac d{(l+2)!}<(2d)^n.
\end{eqnarray}
In the above formula, the factor $(l+2)!$ comes from those factors
$(a_i-z_{\alpha(k+1)})$ for $\alpha=0,1,\cdots,l+1.$ By (\ref{a}) we
know that the  convergence radius of $F(z,t)$ is $>1/(2d)$ for all
$a_i\in A$ as the power series of $t$, which completes the proof.
\end{proof}

\section{$\sigma$\it-convex-sets and Convergence Sets}


{\bD\label{definition} Let  $E$  be a compact subset of $\Omega$.
The holomorphic hull of $E$ in $\Om$ is defined by
$$
\hat{E}_\Om=\{z\in \Omega:|h(z)|\leq\max_{\zeta\in E}|h(\zeta)|\
\hbox{for all}\ h\in \mathscr{O}(\Om)\}.
$$
If $E=\hat{E}_\Om$ we say that $E$ is holomorphically convex in
$\Om$. \eD}

{ \bP \cite[Theorem~1.3.4]{hormander} \label{hull1}A compact subset
$E$ of $\Om\subset \CC$ is holomorphically convex in $\Om$ if and
only if none of the bounded connected components of its complement
is contained in $\Omega$.\eP}
 By Proposition~\ref{hull1} it is easy to prove the following proposition.
 \bP \label{hull2} Let $E$ be a compact subset of $\Om$. Then $\hat{E}_\Om=E\cup
 (\cup_\alpha C_\alpha)$, where $\{C_\alpha\}$ is the set of bounded connected components
 of $\CC\setminus E$ that are contained in $\Om$.\eP

\bL \label{finite} Let $K$ be  a holomorphically convex compact set
in $\Om\subset \CC$. Then $\CC\setminus K$ has only finite components.\eL

\bpp

Suppose that $\CC\setminus K$ has an infinite number of bounded
components. We may  denote by  $C_i,i=1,2,\cdots,$ those components
since $\CC\setminus K$ is an open set. Again denote $d=\max_{w\in
K}|w|$. It is obvious that $\{z\in \CC: |z|>d\}$ must be in the
unbounded component of $\CC\setminus K$. So $C_i\cap\{z\in \CC:
|z|>d\}=\emptyset,\forall i\in 1,2,\cdots $. We get that
 $U=\bigcup_{i=1}^\infty C_i\subset\{z\in
\CC:|z|\leq d\}$. That is, $U$ is bounded.

 By  Proposition
\ref{hull1}, for any $j\in \NN$, $C_j$ can not be contained in
$\Om$.  Then we may choose $w_j\in C_j\setminus \Om$. Let $w$ be a
limited point of the bounded set $\{w_j\}_{j=1}^\infty$. Without the
confusion, we may assume that $w_j\rightarrow w$ when
$j\rightarrow\infty$. Denote by
\begin{equation}\label{inner}
d_j=2{ \sup} \{r>0:  \text{ { There exists $a\in C_j$ such that   }} D(a,r) \subset
C_j \}
\end{equation}
the inner diameter of $C_j$.

By the above argument we get $\sum^\infty_{j=1}{ \frac14}\pi d_j^2\leq
\text{area}(U)\leq \pi d^2$. So we obtain that $d_j$  tends to $0$.
Now by (\ref{inner}), we have
$$
\text{dist}(w_j,K)\leq d_j
$$
because $D(w_j,r)\subset C_j$ for all
$r<\text{dist}(w_j,K)=\text{dist}(w_j,\partial C_j)$. Thus we obtain
that $\text{dist}(w,K)=\lim_{j\rightarrow\infty}\text{dist}(w_j,K)=0$, i.e.,  $w\in K\subset \Omega$.   It is
impossible because $\Omega$  contains only inner points { while} $w$ is a
limited point of { points outer} of $\Om$. Thus $\CC\setminus K$ has
only  finite bounded components.

 \epp


{ \bD A subset $K$ of $\Om\subset \CC$ is said to be a {\it \shc}
set in $\Om$ if it is the union of a countable collection of
holomorphically convex compact subsets of $\Om$. Moreover, $K$ is a
$\sigma$\it-convex-set \rm if it is the union of a countable collection
of polynomial convex compact sets.  \eD}

\bL \label{lunar} Every holomorphically convex compact $K$ in $\Om$
is a $\sigma$-convex-set.\eL

\bpp By Lemma \ref{finite}, $\CC\setminus K$ has only finite
components. We first suppose that it has only one bounded one, named
as $C$. Let $a\in C$ and $ 0<r<R$ such that
$K\subset\mathfrak{R}_{rR}=\{z\in\CC:r<|z-a|<R\}$. For $1<j\in\NN$,
set
$$E_j=\{z=\rho e^{i\theta}+a:\;r<\rho<R;0\leq
\theta\leq \Big(1-\frac{1}{j}\Big)2\pi\}$$
and
$K_j=K\cap E_j$.  Then $K=\cup_{j=2}^\infty K_j$. We need to prove that
$K_j$ is polynomially convex for all $j>1$.

Let $w\notin K_j$ for some one $j\in \NN$. If $w$ is not in $K$, it
should be in $C$ or the unbounded component of $\CC\setminus K${. And  in
either case, } $w$ can be connected with $\infty$ by a path not
intersecting with $K_j$. If $w$ is not in the simply connected set
$E_j$ it is still  path-connected with $\infty$. So $\CC\setminus
K_j$
has no bounded component, that is, $K_j$ is polynomially convex by Proposition \ref{hull1} and
 the remark after Definition \ref{definition}.
\vspace{3mm}

For the case that $K$ is $m$-connected the proof is still valid when we take place
of the fan family $\mathfrak{R}_{rR}=\{z\in\CC:r<|z-a|<R\}\setminus E_j$ by a curved fan family,
each of which cut $K$ into a polynomially convex set.
\epp

 The following proposition  is obvious.

\bP Let $\{K_j\}_{j=1}^\infty$ be a sequence of \shc sets in {
$\Om$}. Then $\cup_{j=1}^\infty K_j$ is
 \shct.\eP

By Lemma \ref{lunar} and the above proposition we have the following
corollary.

 \bC\label{cor} Let $E$ be a  \shc set in { $\Om$}. Then it is also
$\sigma$-convex.\eC

 \vskip2mm
\noindent So we only say the $\sigma$-convex set and not the \shc
set in the rest of this paper.

 \vskip2mm
 Here are two easy examples of  $\sigma$-convex sets. Let
 $\Delta(a,r)=\{z\in\CC:|z-a|<r\}$.
 \bE
Every open set in the complex plane is a $\sigma$-convex  set.
Suppose $E\subset \CC$ is open, then
\begin{equation}
E= \bigcup_{j=1}^\infty \overline{\Delta{(w_j,\frac {r_{w_j}}2)}},
\end{equation}
where $\{w_j\}^\infty_1$ is a dense countable subset of $E$ and
$\frac {r_{w_j}}2=\text{dist}(w_j,\partial E)$.\eE

\bE The unit circle $\Gamma$ is a $\sigma$-convex set  since
$\Gamma=\cup_{j=1}^\infty\{e^{i\theta}:\;0\leq \theta\leq
\Big(1-\frac{1}{j}\Big)2\pi\}$.\eE
The following is a counterexample.
 \bE (Sierpinski
triangle):\rm \ Let { $T$} be the equilateral triangle with
vertices $A,B,C$ and   sides of  length   $1$. And let $D,E,F$ be the mid-points
of sides $AB,BC,AC$ respectively. Denote $4$ equilateral triangles
by{
$$T_{1}=\triangle_{ADF},T_{2}=\triangle_{DBE},T_{3}=\triangle_{ECF},T_{4}=\triangle_{DEF},$$
among which $T_{4}$ }is only the inverted triangle. { Forgot $T_4$ and for $i_1=1,2$ or $3$ we still obtain $T_{i_1}=T_{i_11}\cup T_{i_12}\cup T_{i_12}\cup T_{i_13}\cup T_{i_14}$.}  Continue the process, for  every $T_{i_1\cdots i_k}$ ($i_l=1,2$ or
$3$, $1\leq l\leq k$), we obtain the decomposition: $T_{i_1\cdots
i_k}=T_{i_1\cdots i_k1}\cup T_{i_1\cdots i_k2}\cup T_{i_1\cdots
i_k3}\cup T_{i_1\cdots i_k4}$ where $T_{i_1\cdots i_k4}$ is
inverted one. It is obvious that the side length of $T_{i_1\cdots
i_k}$ is $\frac1{2^{k-1}}$.

The Sierpinski triangle $S$ \rm is defined by
\begin{equation}\label{sier}
  S=\cap_{k=1}^{\infty}\cup_{{ i_l=1,2\text{ or }3}, 1\leq l\leq k}T_{i_1\cdots i_k}=
 {  T\setminus(T_4\cup}\cup_{k=1}^{\infty}\cup_{{ i_l=1,2\text{ or }3}, 1\leq l\leq k}T^o_{i_1\cdots
  i_k4}),
\end{equation}
where $T^o_{i_1\cdots i_k4}$ means the inner part of $T_{i_1\cdots
i_k4}$.\\
 \indent Let $S$ be contained in one
  open set $\Om$ in $\CC$.  The Sierpinski triangle is not a $\sigma$-convex set.  Otherwise, suppose that
 in
$\Om$ there exists a sequence of holomorphically convex compact sets
 $K_j,\ j=1,2,\dots$, such that $S=\cup_{j=1}^\infty K_j$.
But as the closed set of the complete metric space $\Om$, the
Sierpinski triangle is itself a complete metric space,  and hence a
set of second category. By Bair category theorem  there is at least
one $K_j$ containing an induced non-empty open set $V$ in $S$. So there
exists an open subset $U \subset\Om$, such that $V=U\cap S \subset
K_j$. Let $v\in V\subset U$. For every $k\in \NN$ satisfying
$2^{-k}<\text{dist}(v,\partial U)$, by (\ref{sier}), $v$ belongs to
some one $T_{i_1\cdots i_k}$. Then $T_{i_1\cdots i_k}\subset U$. So
$T_{i_1\cdots i_k4}$ is also in $U$. Since $\partial T_{i_1\cdots
i_k4}\subset S$ we get  $\partial T_{i_1\cdots i_k4}\subset V\subset
K_j$. Then $T^o_{i_1\cdots i_k4}$ is one of the complement of $K_j$,
 it is clear that $T^o_{i_1\cdots i_k4}$ is also in $K_j$,
which contradicts with Lemma \ref{finite}. \eE

\bP \label{hull}  Let $K_1$, $K_2$ be compact sets in $\CC$ with
$K_1 \cap K_2=\emptyset$. Then $(K_1\cup K_2) ^\wedge=
\widehat{K}_1\cup \widehat{K}_2$.\eP

\begin{proof}
It suffices to show that $(K_1\cup K_2) ^\wedge\subset
\widehat{K}_1\cup \widehat{K}_2$. Let $U$ be a bounded connected
component of $\CC\setminus ( K_1\cup
 K_2)$, and put $Q=\partial\widehat{\overline{U}}$. Then
  $$
  Q\subset\partial{\overline{U}}\subset\partial U\subset\partial(\CC\setminus(K_1\cup
  K_2))=\partial K_1\cup
  K_2)\subset\overline{K_1\cup K_2}=K_1\cup K_2.
  $$
 Since $Q$ is connected, it follows that $Q\subset K_1$ or $K_2$,
 and hence $U\subset\widehat{Q}\subset \widehat{K}_1\cup \widehat{K}_2$.
 Therefore, $(K_1\cup K_2) ^\wedge\subset
\widehat{K}_1\cup \widehat{K}_2$.
\end{proof}

\bT \label{main} Let $E$ be a $\sigma$-convex  set. Then there exist
polynomially convex compact sets
 $E_n,\;n=1,2,\dots$ in $\Om$, such that $E_{n}\subset E_{n+1}$  for $ n\geq 1$ and $E= \cup_{n=1}^\infty
 E_n$.\eT

\begin{proof}
By the definition of  $\sigma$-convex-set,   $E$ can be written as
\begin{equation}
E =\bigcup_{j=1}^\infty K_{j},
\end{equation} where $K_j$ is the polynomially convex compact set for each
$j$. For $r>0$ we denote by  $N_r(A)$  the $r$-neighborhood of $A$.
For a positive integer $n$ set
 $$F_{n1}= L_{n1}=K_1.$$
 And step by step define
  \bee\label{lll}L_{nj}= K_j\setminus
N_{1/n}\Big(\bigcup_{i=1}^{j-1} K_{i}\Big), \ F_{nj}={\widehat
L}_{nj}, \ \hbox{for}  \ 2\leq j\leq n.\eee Again denote by $$
E_n=\bigcup_{j=1}^n F_{nj}.$$

By Proposition \ref{hull} we obtain
\begin{equation}\label{}
E_n=\bigcup_{j=1}^{n}F_{nj}=\bigcup_{j=1}^{n}{\widehat
L}_{nj}=\widehat{\bigcup_{j=1}^{n} L_{nj}},
\end{equation}
i.e., the set $E_n$ is the polynomial hull of $\bigcup_{j=1}^{n}
L_{nj}$. Since $L_{nj}\subset L_{(n+1)j}$ for each $n\in\NN$ and
$1\leq j\leq n$ we have  $E_n\subset E_{n+1}$.

We now prove that \bee E=\bigcup_{n=1}^{\infty}E_n.\eee Let $z\in
E$. Assume that $j$ is the first positive integer such that $z\in
K_j$, i.e., $z$ in $ K_j$ but not in $K_j,i< j$. Then for
$n>{1}/{\mathrm{dist}(z,\cup_{i=1}^{j-1}K_i)}$ we have  $z\in
K_j\setminus N_{\frac1n}(\cup_{i=1}^{j-1}K_i)=L_{nj}$. And so
$E\subset\bigcup_{n=1}^{\infty}E_n$. On the other hand, let $z\in
\bigcup_{n=1}^{\infty}E_n$. Then there exist $n$ such that $z\in
E_n=\bigcup_{j=1}^{n}{F_{nj}}\subset
\bigcup_{j=1}^{n}\widehat{K_j}=\bigcup_{j=1}^{n}K_j\subset
E$.Therefore $E=\bigcup_{n=1}^{\infty} E_n$.
\end{proof}

 We denote by
$d(E,F)$
 the Euclidean distance between subsets $E$ and $F$.

\bL \label{L314} Let $E, \{K_n\}, \{E_n\}$ be as in Theorem
\ref{main}. Let $U_n= N_{\frac{1}{3n}} (E_n)$. Then for every
positive integer $m$, we have \bee \label{1}{\bigcap_{j=m}^{\infty}
U_j \subset E}.\eee \eL

\begin{proof}
Otherwise, we  suppose that there is an $m>0$ and $z\in\CC$ such
that \bee \label{2} z\in (\bigcap_{j=m}^{\infty} U_j)\setminus
E.\eee

 We first claim that \bee \label{3} z\in N_{\frac{1}{3n}}
\big(\bigcup_{j=1}^{m} F_{nj}\big),\; \hbox{for} \ n=m,
m+1,\dots.\eee We prove (\ref{3}) by induction on $n$. (\ref{3})
obviously holds for $n=m$ by $E_n=\cup_{j=1}^nF_{nj}$ and (\ref{2}).
If there exists a positive integer $N$ so that (\ref{3}) is true for
$n=N-1$ but not $n=N$. Put \bee Q=\bigcup_{j=m+1}^{N} F_{Nj}, \
R=Q\cap \big(\bigcup_{j=1}^{m}F_{Nj}\big), \ \hbox{and} \
S=\big(\bigcup_{j=1}^{m} F_{Nj}\big) \setminus R.\eee

Suppose that $z\in R$. Then there exists $i$ and $k$ with $1\leq
i\leq m<k\leq N$ such that $z\in F_{Ni}\cap F_{Nk}\subset K_i\cap
F_{Nk}$. Hence, $z$ is not in $L_{Nk}$. Since $z\in
F_{Nk}=\widehat{L_{Nk}}$,  we have that $z$ belongs to one of the
bounded components of the complement of $L_{Nk}$. But
$N_{\frac1N}(z)\subset N_{\frac1N}(K_i)$ is contained  in the
complement of $L_{Nk}$. So $N_{\frac1N}(z)$ is just in this bounded
component. Then $N_{\frac1N}(z)$ is contained in the polynomially
convex of $L_{Nk}$, i.e., $N_{\frac1N}(z)\subset
\widehat{L_{Nk}}=F_{Nk}$. So we obtain
\begin{equation}\label{NR}
     N_{\frac{1}{N}}(R)\subset Q.
\end{equation}

Fix $i\leq m$. Let $z$ be in $K_i$ but not in $F_{Nk}$ for $m\leq
k\leq N$. Then $z$ must be in the unbounded component of the
complement of $L_{Nk}$ for $m\leq k\leq N$. Since
$N_{\frac1N}(z)\subset N_{\frac1N}(K_i)$ does not intersect with
$L_{Nk}$, we have that $N_{\frac1N}(z)$ is just in that unbounded
component. That is $N_{\frac1N}(z)\cap \widehat{L_{Nk}}=\emptyset$
for each $m\leq k\leq N$. So we have
\begin{equation}\label{NS}
    d(S,Q)\geq \frac1N.
\end{equation}

Since (\ref{3}) is assumed to hold for $N-1$, and
$N_{\frac{1}{3N-3}} \big(\bigcup_{j=1}^{m} F_{N-1,j}\big)$ is the
subset of
$N_{\frac{1}{3N-3}}(\bigcup_{j=1}^{m}F_{Nj})=N_{\frac{1}{3N-3}}
(S\cup R)$, we see that \bee \label{4} z\in N_{\frac{1}{3N-3}}
(S\cup R).\eee But  by the assumption  that $z\notin E$, (\ref{NR})
and \bee \nonumber N_{\frac{1}{3N-3}} \big( R\big)\subset
N_{\frac{1}{N}}\big(R\big)\subset Q\subset E\eee we obtain
\begin{equation}\label{nnn}
    z\notin N_{\frac{1}{3N}} \big( R\big)\subset N_{\frac{1}{3N-3}} \big(
    R\big).
\end{equation}

 It follows that
\bee \label{5} z\in N_{\frac{1}{3N-3}}(S).\eee By (\ref{NS}) we
obtain
\ba {d(z,Q)&\geq d(Q,S)-d(z,S)\\
&\geq \frac{1}{N}-\frac{1}{3N-3} \geq \frac{1}{3N}\ .} Thus \bee
\label{6} z\notin N_{\frac{1}{3N}} (Q)\ .\eee Combining (\ref{nnn})
with (\ref{6}) contradicts with
 $z\in U_N=N_{\frac{1}{3N}} (Q\cup S)$. Then we complete the induction step for
(\ref{3}). At last, by $\bigcup_{j=1}^{m}F_{nj}\subset
\bigcup_{j=1}^{m}K_j$, we see that \bee \hbox{d}\big(z,
\bigcup_{j=1}^{m} K_j \big) < \frac{1}{3n} \ , \ \forall n\geq m.
\eee Therefore, $z\in \bigcup_{j=1}^{m} K_j\subset E$, contradicting
with (\ref{2}). The proof is complete.
\end{proof}
\bT \label{conv} Let $E$ be a convergence set { in $\Om$}. Then $E$
is a $\sigma$-convex  set. Moreover, there exists an ascending
sequence $\{E_j\}$ of holomorphically convex compact sets in $\Om$
such that $E=\cup_{j=1}^\infty E_j$.\eT

\begin{proof}  Let $E$ be the convergence set $\Conv_\Om(f)$ for
$$
{f(z,t)=f_0(z)+f_1(z)t+\cdots+f_n(z)t^n+\cdots,}$$
 where $f_{n}
(z)\in \mathscr{O}(\Om)$. Denote by
\begin{equation}
E_{jn} = \{z\in\Om:\dist(z,\partial \Om)\ge 1/j, |z|\leq j,
|f_n(z)|\leq j^n\},\;\forall j,n\in\NN,
\end{equation}
where $\dist(z,\partial \Om)$ is the distance between $z$ and the
boundary of $\Om$. Every $E_{jn}$ is obviously compact. Let

\begin{equation}
E_j=\bigcap_{n=1}^\infty E_{jn},\;j\in\NN.
\end{equation}
Then we  have
\begin{equation}
E_j\subset E_{j+1}, \;\;\text{for}\;  j\geq 1.
\end{equation}

By the proof of Theorem~\ref{226}, $E=\cup_{j=1}^\infty E_j$. It
follows from the definition of holomorphically convex sets that each
$E_{jn}$ is holomorphically convex { in $\Om$}. It is also a direct
consequence of the definition that the intersection of a family of
holomorphically convex sets { in $\Om$} is holomorphically convex.
Therefore each $E_j$ is  a holomorphically convex compact set{ in
$\Om$}. The proof of Theorem \ref{conv} is complete.

\end{proof}
For the domain $\Omega$ and a positive integer $m$ denote by
\begin{equation}
\Omega_m = \{z: z\in\Om, \dist(z,\partial \Om)\ge 1/m, |z|\leq m\}.
\end{equation}
\bL \label{9} Let $K$ be a polynomially convex compact subset in
$\Omega$, $U\subset\Omega$ an open set containing $K$, and $m$ a
positive integer. Then there exist a finite number of polynomials, $P_{m1}(z),\dots,P_{m\ell}(z)$, such that
 \begin{equation}\label{last1}
|P_{mj}(z)| \leq 1, \ j=1,\cdots,\ell, \ \mathrm{for\ all}\ z\in K,
\end{equation}
 and
\begin{equation}\label{last2}
\max_{j}\{|P_{m1}(z)|,\dots,|P_{m\ell}(z)|\}\geq m \ , \ \mathrm{for\ all}\
z\in \Omega_m\setminus U.
 \end{equation} \eL
\begin{proof}

Due to the polynomial convexity of $K$, for each $z_0\in \Om_m \setminus U$ there
exists a polynomial $Q(z)$ such that
$$
|Q(z)|\leq 1, \ \mathrm{for} \ z\in K, \ \ \hbox{and} \ \ |Q(z_0)|\geq m+1.
$$
 Then there is some
neighborhood $V(z_0) \ \hbox{of}\ z_0$ such that $  |Q(z)|\geq m  \
\hbox{for each} \ z\in V(z_0)$. Since $\Om_m \setminus U$ is compact there are a
finite number of such open sets $V(z_1),\dots,V(z_{\ell})$ covering
$\Om_m \setminus U$. The corresponding polynomials are denoted by  $P_{mj}(z), \
\forall j=1,\dots, \ell$. Then $P_{mj}(z), \
\forall j=1,\dots, \ell$, satisfy (\ref{last1}) and (\ref{last2}).
\end{proof}

\rm Now we prove the main theorem of this paper.

\bT \label{quas} $E\subset\Om$ is a convergence set { in $\Om$} if
and only if it is
 $\sigma$-convex.\eT

\begin{proof} By Corollary~ \ref{cor} and Theorem~ \ref{conv} we only need to prove 'if'.\\
\indent Let $E$ be a $\sigma$-convex  set. By Theorem~\ref{main},
there exist polynomially convex compact sets $E_n$ such that
$E=\bigcup_{n=1}^{\infty} E_n$, and $E_n\subset E_{n+1} \ \hbox{for}
\ n\geq 1$. Let $U_n$ be the neighborhood of $E_n$  in Theorem~\ref{L314}, then  we have
$$\bigcup_{k=1}^{\infty}\bigcap_{n=k}^{\infty} U_n\subset E. $$
And by $E_k\subset E_{k+1}$ and $E_k\subset U_n$ for $n\geq k$ we
obtain
 \bee \label{22}
 \nonumber E=\bigcup_{k=1}^{\infty}\bigcap_{n=k}^{\infty} U_n.\eee

Now for each $k\in \NN$, let $P_{k1},\cdots,P_{kn_k}$ be the
polynomials for $\Om_m=\Om_k,K=E_k$ and $U=U_k$ by Lemma
\ref{9}. Then we have
$$
|P_{ks}(z)| \leq 1 \  ,  \mathrm{for} \ z\in E_k,\  1\leq s\leq n_k
$$
 and for every
$z\in \Omega_{k}\setminus U_k$ there exists  $j$, $1\leq j\leq
n_k$, such that
 $$ |P_{kj}(z)|\ge k.$$
  Enumerate the countable set
of polynomials
$\{\{P_{kj}\}_{j=1}^{n_k}\}_{k=1}^{\infty}$ by $\{h_\ell\}_{\ell=1}^{\infty}$
such that
$h_1=P_{11},h_2=P_{12},\cdots,h_{n_1}=P_{1n_1},h_{n_1+1}=P_{21},\cdots$.
Define
$$
 f_0(z)=1$$
 and
 $$
  f_{\ell}(z)=h_{\ell}^{\ell}(z) \ , \ \forall
\ell\in\NN^+.
$$

For any $z\in E$  there exist $k$ such that $z\in E_n$ for $n\geq
k$. Then for $\ell>n_1+\cdots+n_k$ we have  $|f_\ell(z)|\leq
|h_\ell(z)| \leq1$.  {It implies that $z\in \Conv_\Om(f)$. Thus
$E\subset \Conv_\Om(f)$.}

For any $z\in \Om\setminus
E=\Om\setminus(\bigcup_{k=1}^{\infty}\bigcap_{n=k}^{\infty} U_n)=\bigcap_{k=1}^{\infty}\bigcup_{n=k}^{\infty} (\Om\setminus U_n)=\bigcap_{n=1}^{\infty}(\Om\setminus U_n)$, chose any positive integer
$m$ such that $z\in \Om_{m}$. Then for any $k>m$ we have $z\in \Om_{k}$. By Lemma \ref{9}  and
taking $\Om_m=\Om_{k},K=E_k, U=U_k$, we have a polynomial
$P_{kj}(z)$ such that $|P_{kj}(z)|> k$. Let
$\ell=n_1+\cdots+n_{k-1}+j$.  We obtain
$$
 |f_l(z)|> k^l\geq
m^l.
$$
 To summarize, for each $z\notin E$, and for each
positive integer $m$, there exist infinite positive integers  $\ell$ such that
$|f_{\ell}(z)|>m^{\ell}$. Consequently, it implies that { $z\notin
\Conv_\Om(f)$. Hence $\Conv_\Om(f)\subset E$. Therefore
$E=\Conv_\Om(f)$.}
\end{proof}

\noindent It is difficult to directly deal with the union of a
countable collection of convergence sets { in $\Om$}. But by
Theorem~\ref{quas} it is equivalent to the case of $\sigma$-convex
sets, while the latter is easily done by definition.
 \bC The union of a
countable collection of convergence sets { in $\Om$} is a
convergence set { in $\Om$}.\eC
\brs
  In this paper we don't confine the degree of the coefficient polynomials in the definition of the convergence set.  However  it is mentioned in the first section that $\Conv(f)$ is a polar set if the degree of coefficients $P_n(z)$  are assumed  not great than $n$. While it is surprising that it is not essential for  the degree of the coefficient polynomial when it is great than $n$.
   For $\varepsilon>0$, let the $n$-th term coefficients  of   $f(z,t)$ is polynomial with degree not great than $n^{1+\varepsilon}$. We call $\Conv(f)$ as the $\varepsilon$-convergence-set.
Then the $E\subset \CC$ is a convergence set if only if it is an $\varepsilon$-convergence-set. In fact,
Let $E=\Conv(f)$, where$
  f(z,t)=\sum_{j=0}^\infty P_n(z)t^n$.
Denote by $d_n$ the degree of $P_n$, $n=0,1,2,\cdots$. We take  $m_{n}$ be a sequence of increasing integers satisfy that $m_n>(\frac{d_n}{n^{1+\varepsilon}})^{\frac1\varepsilon}$. Now let
$F(z,t)=\sum_{j=0}^\infty (P_n(z))^{m_n}t^{nm_n}.
$
Then  $E=\Conv(f)=\Conv(F)$, and $E=\Conv(F)$ be an $\varepsilon$-convergence-set.
\ers
{\bf Acknowledgment.} We thank Buma Fridman for helpful discussions. Part of the third named author's work was done while visiting Tshinghua University Yau Mathematical Sciences Center during his sabbatical leave
in spring 2014. He is grateful for the Center's hospitality and financial support.

\end{document}